# Modified parameter of the Dai–Liao conjugacy condition of the conjugate gradient method with some applications

[1*]Ahmad Alhawarat, [2]Sultanah Masmali, [1]Ali Jaradat, [2]Ramadan Sabra, [3]Shahrina Ismail

[1]Department of Mathematics, Faculty of Arts and Science, Amman Arab University, Amman 11953, Jordan
[2]Department of Mathematics, College of Science, Jazan University, Jazan, Saudi Arabia.
[3]Financial Mathematics Program, Faculty of Science and Technology, Universiti Sains Islam Malaysia, Bandar Baru Nilai, 71800, Nilai, Negeri Sembilan, Malaysia.

\* **Correspondence:** Email:abadee2010@yahoo.com

**Abstract**

In this research, we propose a modified parameter of the Dai–Liao conjugacy condition of the CG method with the restart property. The proposed method depends on the Lipschitz constant and is related to the Hestenes–Stiefel method and satisfies the descent condition and global convergence properties for convex and non-convex functions. In performed numerical experiments, we compare the new method with CG Descent using more than 200 functions from the CUTEst library. The comparison results show that the new method outperforms CG Descent in terms of CPU time, number of iterations, number of gradient evaluations, and number of function evaluations. Finally, we present some applications on heat conduction problem and image restoration.

## 1- Introduction

Conjugate gradient (CG) methods have been widely used for solving nonlinear unconstrained optimization problems due to their low memory requirements for implementation. Moreover, CG methods have been used in many applications such as regression analysis, image restoration, electrical circuits, and many others. CG method is used to determine optimal solutions for the following optimization problem

$$\min f(x),\ x \in \mathbb{R}^n,$$

where $f: \mathbb{R}^n \to \mathbb{R}$ is a continuously differentiable function, and its gradient $\nabla f(x_k) = g_k = g(x_k)$ should exist. From the starting point (arbitrary or standard) $x_1 \in \mathbb{R}^n$, the CG method generates a sequence of vectors $x_k$ as follows:

$$x_{k+1} = x_k + \alpha_k d_k, k = 1, 2, \ldots, \quad (1)$$

where $x_k$ represents the present iteration and $\alpha_k > 0$ represents a step size obtained from a line search, such as the exact line search or inexact line search, which will be discussed in the succeeding sections. The search direction $d_k$ of the CG method is defined by:

$$d_k = \begin{cases} -g_k, & \text{if } k = 1, \\ -g_k + \beta_k d_{k-1}, & \text{if } k \geq 2, \end{cases} \quad (2)$$

where $\beta_k$ is the update parameter. The following exact line search can be utilized to obtain the step size $\alpha_k$,

$$f(x_k + \alpha_k d_k) = \min f(x_k + \alpha d_k), \alpha \geq 0. \qquad (3)$$

However, Eq. (3) is computationally expensive because it requires unidimensional optimization to achieve the step size and many iterations to reach convergence. To avoid this problem, the inexact line search is dominant approach in computing the step size. The most popular inexact line search is the strong Wolfe–Powell (SWP) line search [1, 2], which is defined as

$$f(x_k + \alpha_k d_k) \leq f(x_k) + \delta \alpha_k g_k^T d_k \qquad (4)$$

and

$$|g(x_k + \alpha_k d_k)^T d_k| \leq \sigma |g_k^T d_k|, \qquad (5)$$

where $0 < \delta < \frac{1}{2}, \delta < \sigma < 1$

A version of the Wolfe–Powell line search is the weak Wolfe–Powell (WWP) line search, which is defined by (4) and

$$g(x_k + \alpha_k d_k)^T d_k \geq \sigma g_k^T d_k.$$

The most famous classical formulae of the CG methods are the Hestenes–Stiefel (HS)[3], Fletcher–Reeves(FR)[4], and Polak–Ribiere–Polyak (PRP) methods [5], which are defined by parameters:

$$\beta_k^{HS} = \frac{g_k^T y_{k-1}}{d_{k-1}^T y_{k-1}}, \quad \beta_k^{FR} = \frac{\|g_k\|^2}{\|g_{k-1}\|^2} \text{ and } \beta_k^{PRP} = \frac{g_k^T y_{k-1}}{\|g_{k-1}\|^2}, \text{ where } y_{k-1} = g_k - g_{k-1}.$$

Powell in [6] provided a counterexample showing that a non-convex function exists in which the PRP and HS methods fail to satisfy the convergence properties, even if the exact line search is employed. Powell recommended using the nonnegative of $\beta_k^{HS}$ and $\beta_k^{PRP}$ parameters to obtain the convergence properties of CG method. Gilbert and Nocedal [7] proved that the nonnegative PRP or HS method defined by $\beta_k = \max\{\beta_k^{PRP}, 0\}$, is globally convergent with arbitrary line searches.

The descent condition (downhill condition) plays a crucial role in the convergence of the CG method and its robustness, and it is defined as follows:

$$g_k^T d_k < 0. \qquad (6)$$

Al-Baali [8] proposed another version of the downhill condition called the sufficient descent condition, which also plays a significant role in the convergence of the CG method. Al-Baali used Eq. (7) to establish the global convergence properties of $\beta_k^{FR}$, which is defined as follows:

If there exists a constant $c > 0$, $g_k^T d_k \leq -c\|g_k\|^2$, $\forall k \in N$, then the search direction $d_k$ guarantees the sufficient descent condition.

$$g_k^T d_k \leq -c\|g_k\|^2 \qquad (7)$$

Based on the quasi-Newton method, the Broyden–Fletcher–Goldfarb–Shanno (BFGS) method and the limited-memory BFGS (LBFGS) method, and using Eq. (2), Dai and Liao[9] proposed the following conjugacy condition:

$$d_k^T y_{k-1} = -t g_k^T s_{k-1}, \qquad (8)$$

where $s_{k-1} = x_k - x_{k-1}$, and $t \geq 0$. In the case of $t = 0$, Eq. (8) is considered the classical conjugacy condition. Using Eqs. (2) and (8), Dai and Liao[9] proposed the following CG formula:

$$\beta_k^{DL} = \frac{g_k^T y_{k-1}}{d_{k-1}^T y_{k-1}} - t \frac{g_k^T s_{k-1}}{d_{k-1}^T y_{k-1}}. \qquad (9)$$

However, $\beta_k^{DL}$ can not satisfy the descent condition and convergence properties similar to $\beta_k^{PRP}$ and $\beta_k^{HS}$, i.e., $\beta_k^{DL}$ is not nonnegative in general. Thus, Dai and Liao [9] replaced Eq. (9) with the following equation:

$$\beta_k^{DL+} = \max\{\beta_k^{HS}, 0\} - t \frac{g_k^T s_{k-1}}{d_{k-1}^T y_{k-1}}. \tag{10}$$

However, $\beta_k^{DL+}$ cannot satisfy the descent property in some cases. Therefore, Dai and Liao[9] restarted (10) using a negative gradient(steepest descent) when $\beta_k^{DL+}$ fails to satisfy (7). Another method for determining the optimal parameter $t$ was proposed by Babaie-Kafaki and Ghanbari [10, 11], where they rewrote the search direction (Eq. (2)) with $\beta_k^{DL}$ and based on Perry [12], as follows: $d_{k+1} = -Q_{k+1} g_{k+1}$, where $Q_{k+1} = I - \frac{s_k y_k^T}{s_k^T y_k} + t \frac{s_k s_k^T}{s_k^T y_k}$. Babaie-Kafaki and Ghanbari[10] proposed the following adaptive choices for $t$:

$$t = \frac{s_k y_k^T}{\|s_k\|^2} + \frac{\|y_k\|}{\|s_k\|}, \text{ and } t = \frac{\|y_k\|}{\|s_k\|}.$$

Andrei in [13] proposed CG method with the following parameter:

$$\beta_k^{DL*} = \max\left\{\frac{y_k^T g_k}{y_k^T s_k}, 0\right\} - t_k^* \frac{s_k^T g_{k+1}}{y_k^T s_k},$$

where $t_k^* = y_k^T s_k / \|s_k\|^2$. Hager and Zhang [14,15] presented a modified CG parameter that satisfies the descent property for any inexact line search with $g_k^T d_k \leq -(7/8)\|g_k\|^2$. This new version of the CG method is globally convergent whenever the line search satisfies the WWP line search. This formula is expressed as follows:

$$\beta_k^{HZ} = \max\{\beta_k^N, \eta_k\},$$

where $\beta_k^N = \frac{1}{d_k^T y_k}(y_k - 2d_k \frac{\|y_k\|^2}{d_k^T y_k})^T g_k$, $\eta_k = -\frac{1}{\|d_k\| \min\{\eta, \|g_k\|\}}$, and $\eta > 0$ is a constant.

Note that, if $t = 2\frac{\|y_k\|^2}{s_k^T y_k}$, then $\beta_k^N = \beta_k^{DL}$. Zhang et al. [16] proposed a new parameter for Eq. (9) as follows:

$$t = \frac{\|y_k\|^2}{s_k^T y_k} - \frac{1}{4}\frac{s_k^T y_k}{\|s_k\|^2}.$$

Yao et al. [17] proposed three terms of CG with a new choice of $t$ as follows:

$$d_{k+1} = -g_{k+1} + \left(\frac{g_k^T y_k - t_k g_{k+1}^T s_k}{y_k^T d_k}\right) d_k + \frac{g_{k+1}^T d_k}{y_k^T d_k} y_k.$$

Based on the SWP line search, Yao *et al.* [17] selected $t_k$ to satisfy the descent condition

$$t_k > \frac{\|y_k\|^2}{y_k^T s_k}.$$

Yao et al. [17] also proposed a theorem stating that if $t_k$ is close to $\frac{\|y_k\|^2}{y_k^T s_k}$, then the search direction results in a zigzag search path. Therefore, they selected the following choice for $t_k$:

$$t_k = 1 + 2\frac{\|y_k\|^2}{y_k^T s_k}.$$

Al-Baali et al. [18] proposed a new CG version called G3TCG that offers many selections of CG parameters. They found that the G3TCG method is more efficient than $\beta_k^{HZ}$ in some cases and competitive in some other cases.

There are many applications of CG methods in several fields such as electrical engineering, image restorations, machine learning, economics, and many others. We advise the reader to see the following references [19-23].

## 2- Proposed CG formula and its motivation

The CG method with $\beta_k^{DL}$ cannot satisfy the descent condition, but $\beta_k^{DL}$ inherits the conjugacy condition. To improve the properties of $\beta_k^{DL}$, we used $\beta_k^{AZPRP}$ as presented by Alhawarat et al. [24] to propose a new nonnegative CG method that can satisfy the sufficient descent condition and global convergence properties with the SWP line search as follows:

$$\beta_k^{AZPRP} = \begin{cases} \dfrac{\|g_k\|^2 - \mu_k|g_k^T g_{k-1}|}{\|g_{k-1}\|^2}, & \|g_k\|^2 > \mu_k|g_k^T g_{k-1}| \\ 0, & \text{otherwise.} \end{cases}$$

The new formula is a modification of $\boldsymbol{\beta_k^{DL}}$ and $\boldsymbol{\beta_k^{HS}}$, with the restart criterion depending on the Lipschitz constant used in the study conducted by Alhawarat et al. [24]. The modified formula is expressed as follows:

$$\beta_k^{AZHS} = \begin{cases} \dfrac{\|g_k\|^2 - \mu_k|g_k^T g_{k-1}|}{d_{k-1}^T y_{k-1}} - \dfrac{1}{\alpha_k}\mu_k \dfrac{g_k^T s_{k-1}}{d_{k-1}^T y_{k-1}}, & \text{if } \|g_k\|^2 > \mu_k|g_k^T g_{k-1}|, \\ -\dfrac{1}{\alpha_k}\mu_k \dfrac{g_k^T s_{k-1}}{d_{k-1}^T y_{k-1}}, & \text{otherwise,} \end{cases} \quad (11)$$

where $\|\cdot\|$ represents the Euclidean norm and $\mu_k$ is defined as follows:

$$\mu_k = \dfrac{\|s_{k-1}\|}{\|y_{k-1}\|}.$$

In the first case of equation (11) we can note that

$$\beta_k^{AZHS} \leq \dfrac{\|g_k\|^2}{d_{k-1}^T y_{k-1}} - \dfrac{1}{\alpha_k}\mu_k \dfrac{g_k^T s_{k-1}}{d_{k-1}^T y_{k-1}} \quad (12)$$

It is worth noting that the assignment (11) inherits the advantages of $\beta_k^{DL}$, $\beta_k^{HS}$ and $\beta_k^{AZPRP}$. Usage of the proposed parameter $\beta_k^{AZHS}$ in (11) leads to the novel CG method described in Algorithm 1.

**Algorithm 1**
**Step 1**. Set a starting point $x_1$, the initial point can be arbitrary or standard for scientific functions. The initial search direction is the negative gradient, i.e., $d_1 = -g_1$. Let $k := 1$
**Step 2**. If the stopping condition is satisfied, then stop.
**Step 3**. Compute the search direction $d_k$ based on Eq. (2) using Eq. (11).
**Step 4**. Compute the step size $\alpha_k$ using Eqs. (4) and (5).
**Step 5.** Update $x_{k+1}$ based on Eq. (1).
**Step 6**. Set $k := k + 1$ and go to Step 2.

### 3- Convergence analysis of $\beta_k^{AZHS}$

To satisfy the convergence analysis of the modified CG method, we consider the following assumptions:

**Assumption 1**
**A**. The level set $\Phi = \{x|f(x) \leq f(x_1)\}$ is bounded. In other words, a positive constant $B$ exists, such that
$$\|x\| \leq B, \forall x \in \Phi.$$
**B**. In some neighbourhood $P$ of $\Phi$, $f$ is continuously differentiable, and its gradient is Lipschitz continuous. In other words, for all $x, y \in P$, there exists a constant $L > 0$, such that
$$\|g(x) - g(y)\| \leq L\|x - y\|.$$

This assumption implies that there exists a positive constant $\bar{\gamma}$, such that
$$\|g(x)\| \leq \bar{\gamma}, \quad \forall x \in N.$$

**Theorem 3.1** Let the sequences $\{g_k\}$ and $\{d_k\}$ be obtained using Eqs. (1) and (2), and $\beta_k^{AZHS}$ where $\alpha_k$ is computed using the SWP line search in Eqs. (4) and (5). If $\sigma \in (0, 0.5)$, then the descent condition provided in Eq. (7) holds.

**Proof.** The proof recognizes two cases.
**Case 1:** $\|g_k\|^2 > \mu_k |g_k^T g_{k-1}|$. This assumption implies
$$\beta_k^{AZHS} = \frac{\|g_k\|^2 - \mu_k |g_k^T g_{k-1}|}{d_{k-1}^T y_{k-1}} - \frac{1}{\alpha_k} \mu_k \frac{g_k^T s_{k-1}}{d_{k-1}^T y_{k-1}}.$$
Multiplying Eq. (2) with $g_k^T$ it can be concluded
$$g_k^T d_k = g_k^T(-g_k + \beta_k d_{k-1}) = -\|g_k\|^2 + \beta_k g_k^T d_{k-1}.$$

$$\leq -\|g_k\|^2 + \frac{\|g_k\|^2}{|d_{k-1}^T y_{k-1}|} |g_k^T d_{k-1}| - \mu_k \frac{\|g_k^T d_{k-1}\|^2}{d_{k-1}^T y_{k-1}}.$$

Using the SWP line search we obtain the following inequality
$$\frac{|g_k^T d_{k-1}|}{|d_{k-1}^T y_{k-1}|} \leq \frac{\sigma}{1-\sigma}.$$

Thus,
$$g_k^T d_k \leq -\|g_k\|^2 + \frac{\sigma \|g_k\|^2}{(1-\sigma)}.$$
$$\leq -\|g_k\|^2 + \frac{\sigma \|g_k\|^2}{(1-\sigma)} = -\|g_k\|^2 \left(1 - \frac{\sigma}{1-\sigma}\right).$$

Let $c = \left(1 - \frac{\sigma}{1-\sigma}\right)$. In this case, if $\sigma < \frac{1}{2}$ one concludes
$$g_k^T d_k \leq -c \|g_k\|^2.$$

**Case 2:** $\|g_k\|^2 \leq \mu_k |g_k^T g_{k-1}|$
This assumption implies
$$\beta_k^{AZHS} = -\frac{1}{\alpha_k} \mu_k \frac{g_k^T s_{k-1}}{d_{k-1}^T y_{k-1}}$$
and further
$$g_k^T d_k = g_k^T(-g_k + \beta_k d_{k-1}) = -\|g_k\|^2 + \beta_k g_k^T d_{k-1}.$$

$$\leq -\|g_k\|^2 + \left(-\frac{\mu_k}{\alpha_{k-1}} \frac{g_k^T s_{k-1}}{d_{k-1}^T y_{k-1}}\right) g_k^T d_{k-1}.$$
$$= -\|g_k\|^2 - \mu_k \frac{\|g_k^T d_{k-1}\|^2}{d_{k-1}^T y_{k-1}}.$$

Since we use SWP line search, we obtain $d_{k-1}^T y_{k-1} > 0$, thus
$$g_k^T d_k \leq -c \|g_k\|^2.$$
which finishes the proof.

The following Lemma shows that if $L > 1$ then equation (13) holds. Note that if $L << 1$ then $\|g_k\|^2 > \mu_k |g_k^T g_{k-1}|$ can not be satisfied.

**Lemma 3.1** If $\|g_k\|^2 > \mu_k |g_k^T g_{k-1}|$ and $L > 1$, then
$$\|g_k\|^2 - \frac{1}{L} |g_k^T g_{k-1}| \leq L \left| \|g_k\|^2 - |g_k^T g_{k-1}| \right|. \tag{13}$$

**Proof.** We will perform the proof using contradiction. Suppose that

$$\|g_k\|^2 - \frac{1}{L}|g_k^T g_{k-1}| > L\left|\|g_k\|^2 - |g_k^T g_{k-1}|\right|,$$

and divide both sides by $L$:

$$\frac{\|g_k\|^2}{L} - \frac{1}{L^2}|g_k^T g_{k-1}| > \left|\|g_k\|^2 - |g_k^T g_{k-1}|\right| \qquad (14)$$

Using Assumption 1, we derive the following:

$$\|g_k\|^2 > \mu_k |g_k^T g_{k-1}| > \frac{1}{L}|g_k^T g_{k-1}|.$$

If $L > 1$, we conclude that inequality (14) is not true, and this results into contradiction. Thus, inequality (13) is true.

The following Lemma 3.2 indicates the step length always has a lower bound.

**Lemma 3.2**[25]. Suppose that the objective function satisfies Assumption 1. If the step length $\alpha_k$ fulfils the SWP line search conditions (4) and (5), then

$$\alpha_k \geq \frac{(1-\sigma)|g_k^T d_k|}{L\|d_k\|^2}.$$

The condition presented in Eq. (15) is called the Zoutendijk condition [26] and plays an important role in proving the convergence properties of the CG method. We use the contradiction technique with Eq. (15) to prove that $\liminf_{k\to\infty}\|g_k\|=0$.

**Lemma 3.3** Assume that Assumption 1 holds. Consider any form of Eqs. (1) and (2), with step size $\alpha_k$ satisfying the WWP line search, where the search direction $d_k$ is descent. We can obtain the following:

$$\sum_{k=0}^{\infty} \frac{(g_k^T d_k)^2}{\|d_k\|^2} < \infty. \qquad (15)$$

Moreover, Eq. (15) holds for the exact and SWP line searches. By substituting Eq. (7) into Eq. (15), we obtain the following:

$$\sum_{k=0}^{\infty} \frac{\|g_k\|^4}{\|d_k\|^2} < \infty. \qquad (16)$$

**Lemma 3.4** if $\|g_k\|^2 > \mu_k |g_k^T g_{k-1}|$ satisfied then $\mu_k = \frac{\|s_{k-1}\|}{\|y_{k-1}\|}$ is bonded above and below.

**Proof.** Since $\|g_k\|^2 > \mu_k |g_k^T g_{k-1}| > \frac{1}{L}|g_k^T g_{k-1}|$, based on Assumption 1 we conclude $0 < \mu_k \leq E$, where $E$ denotes a positive constant. Moreover, if $y_{k+1} = 0$, this means $x_{k+1} = x_k$ and we know that $x_{k+1} = x_k + \alpha_k d_k$, Thus $\alpha_k d_k = 0$. However, by Lemma 3.2 we conclude that $\alpha_k > 0$. This mean $d_k = 0$. By using Theorem 3.1 and Lemma 3.3 we have a contradiction.

Dai *et al.* [27] presented the following theorem, which is also useful for proving the global convergence properties of CG methods.

**Theorem 3.2.** Suppose that Assumption 1 holds. Consider any CG method in the forms of Eqs. (1) and (2), where $d_k$ is a descent direction and $\alpha_k$ is obtained using the SWP line search. If

$$\sum_{k\geq 1}^{\infty} \frac{1}{\|d_k\|^2} = \infty,$$

then

$$\liminf_{k\to\infty}\|g_k\|=0.$$

### *3.1 Global convergence properties for the convex function*

In the following theorem, if $f(x)$ is a uniformly convex function, then the CG method satisfies $\beta_k^{AZHS}$ strong global convergence properties.

**Theorem 3.3.** Suppose that Assumption1 holds. Consider the CG method in the forms of Eqs. (1) and (2), with $\beta_k^{AZHS}$, $L > 1$, and $d_k$ as a descent direction, where $\alpha_k$ is obtained using the SWP line search. If $f(x)$ is a uniformly convex function, then $\liminf_{k \to \infty} \|g_k\| = 0$.

**Proof.** Because the function $f(x)$ is uniformly convex, there exists a positive constant $\varpi$:
$$\varpi \|x - y\|^2 \leq (\nabla f(x) - \nabla f(y))^T (x - y).$$
For all $x, y \in P$. Thus,
$$d_{k-1} y_{k-1} \geq \varpi \alpha_{k-1} \|d_{k-1}\|^2 \tag{17}$$

$$\beta_k^{AZHS} = \frac{\|g_k\|^2 - \mu_k |g_k^T g_{k-1}|}{d_{k-1}^T y_{k-1}} - \frac{\mu_k}{\alpha_{k-1}} \frac{g_k^T s_{k-1}}{d_{k-1}^T y_{k-1}}$$
$$\leq \frac{\|g_k\|^2 - \frac{1}{L}|g_k^T g_{k-1}|}{d_{k-1}^T y_{k-1}} - \frac{\mu_k}{\alpha_{k-1}} \frac{g_k^T s_{k-1}}{d_{k-1}^T y_{k-1}}$$

Using Eqs. (17) and (13), we obtain
$$\beta_k^{AZHS} \leq \frac{L\|g_k\|(\|g_k - g_{k-1}\|)}{\varpi \alpha_{k-1} \|d_{k-1}\|^2} + E \frac{|g_k^T s_{k-1}|}{\varpi \alpha_{k-1}^2 \|d_{k-1}\|^2}$$
$$\leq \frac{L\|g_k\|\|g_k - g_{k-1}\|}{\varpi \alpha_{k-1} \|d_{k-1}\|^2} + E \frac{\|g_k\|\|s_{k-1}\|}{\varpi \alpha_{k-1}^2 \|d_{k-1}\|^2}$$

Applying Assumption 1, we obtain the following:
$$\beta_k^{AZHS} \leq \frac{L^2 \|g_k\| \alpha_{k-1} \|d_{k-1}\|}{\varpi \alpha_{k-1} \|d_{k-1}\|^2} + E \frac{\|g_k\|\|d_{k-1}\|}{\varpi \alpha_{k-1} \|d_{k-1}\|^2}$$
$$\leq \frac{L^2 \|g_k\|}{\varpi \|d_{k-1}\|} + E \frac{\|g_k\|}{\varpi \alpha_{k-1} \|d_{k-1}\|} = \frac{\|g_k\|}{\|d_{k-1}\|} \left( \frac{L^2}{\varpi} + \frac{E}{\varpi \alpha_{k-1}} \right)$$

Based on Eq. (2), we obtain
$$\|d_k\| \leq \|g_k\| + |\beta_k|\|d_{k-1}\|$$
$$\leq \|g_k\| + \frac{\|g_k\|}{\|d_{k-1}\|} \left( \frac{L^2}{\varpi} + \frac{E}{\varpi \alpha_{k-1}} \right) \|d_{k-1}\|$$
$$\leq \bar{\gamma} \left( 1 + \left( \frac{L^2}{\varpi} + \frac{E}{\varpi \alpha_{k-1}} \right) \right).$$

Thus, by using Theorem 3.2, we obtain the following:
$$\liminf_{k \to \infty} \|g_k\| = 0.$$

### 3.2 Global convergence for $\beta_k^{AZHS}$ with the SWP line search for general functions

Using Property(*) and some lemmas, Gilbert and Nocedal [7] proved the global convergence of nonnegative PRP and HS methods. Because our modification is nonnegative and satisfies Property(*), by using the other lemmas presented below, we perform our proof, like [7]. This property is defined as follows:

**Property(*)**
Consider any CG method in the form of Eqs. (1) and (2). Assume
$$0 < \gamma \leq \|g_k\| \leq \bar{\gamma} \tag{18}$$

for all $k \geq 1$. The CG method then inherits Property(*) if for $\forall k$, there exist constants $b > 1$ and $\lambda > 0$, such that $|\beta_k| \leq b$ and $\|s_k\| \leq \lambda$, which implies that $|\beta_k| \leq \frac{1}{2b}$.
The following lemma shows that $\beta_k^{AZHS}$ satisfies Property(*).

**Lemma 3.5** Consider a CG method in the form of Eqs. (1) and (2) using $\beta_k^{AZHS}$ with $L > 1$; and Lemma 3.1 holds true then, $\beta_k^{AZHS}$ satisfies Property(*).

**Proof.** Let $b = \frac{2L\alpha_{k-1}\bar{\gamma}^2 + B\bar{\gamma}}{\alpha_{k-1}L(1-\sigma)c\gamma^2} \geq 1$, and let $\lambda \leq \frac{(1-\sigma)c\gamma^2}{2(L^2 + \frac{E}{\alpha_{k-1}})\bar{\gamma}b}$. Then the following inequality holds:

$$\beta_k^{AZHS} \leq \frac{\|g_k\|^2 - \mu_k |g_k^T g_{k-1}|}{d_{k-1}^T y_{k-1}} - \frac{\mu_k}{\alpha_{k-1}} \frac{g_k^T s_{k-1}}{d_{k-1}^T y_{k-1}}.$$

Using Eqs. (13) and (18), we obtain the following:

$$\beta_k^{AZHS} \leq \frac{\|g_k\|^2 + |g_k^T g_{k-1}|}{d_{k-1}^T y_{k-1}} + \frac{E}{\alpha_{k-1}} \frac{|g_k^T s_{k-1}|}{d_{k-1}^T y_{k-1}}$$
$$\leq \frac{2\bar{\gamma}^2}{(1-\sigma)c\gamma^2} + \frac{EB\bar{\gamma}}{\alpha_{k-1}L(1-\sigma)c\gamma^2} = \frac{2L\alpha_{k-1}\bar{\gamma}^2 + EB\bar{\gamma}}{\alpha_{k-1}(1-\sigma)c\gamma^2} = b > 1.$$

If $\|s_k\| \leq \lambda$, we obtain the following:

$$\beta_k^{AZHS} \leq \frac{\|g_k\|^2 - \mu_k |g_k^T g_{k-1}|}{d_{k-1}^T y_{k-1}} - \frac{\mu_k}{\alpha_{k-1}} \frac{g_k^T s_{k-1}}{d_{k-1}^T y_{k-1}} \leq \frac{L\|g_k\| \|(g_k - g_{k-1})\|}{d_{k-1}^T y_{k-1}} + \frac{E}{\alpha_{k-1}} \frac{\|g_k\| \|s_{k-1}\|}{d_{k-1}^T y_{k-1}}$$
$$\leq \frac{L^2\|g_k\| \|s_{k-1}\|}{d_{k-1}^T y_{k-1}} + \frac{E}{\alpha_{k-1}} \frac{\|g_k\| \|s_{k-1}\|}{d_{k-1}^T y_{k-1}}$$
$$\leq \frac{(L^2 + \frac{E}{\alpha_{k-1}})\|g_k\| \|s_{k-1}\|}{d_{k-1}^T y_{k-1}} \leq \frac{(L^2 + \frac{E}{\alpha_{k-1}})\bar{\gamma}\lambda}{(1-\sigma)c\gamma^2} = \frac{1}{2b}.$$

Thus, the proof is complete.

The following lemmas match with Lemma 4.1 and Lemma 4.2, as presented by Gilbert and Nocedal[7].

**Lemma 3.6** Assume that Assumption 1 holds and the sequences $\{g_k\}$ and $\{d_k\}$ are generated using Algorithm 1, where the step size $\alpha_k$ is computed via the SWP line search, such that the sufficient descent condition holds. If $\beta_k \geq 0$, there exists a constant $\gamma > 0$, such that $\|g_k\| > \gamma$ for all $k \geq 1$. Then, $d_k \neq 0$ and

$$\sum_{k=0}^{\infty} \|u_{k+1} - u_k\|^2 < \infty, \text{ where } u_k = \frac{d_k}{\|d_k\|}.$$

**Proof.** First, if $d_k = 0$, then from the sufficient descent condition, we obtain $g_k = 0$. Thus, we suppose $d_k \neq 0$ and

$$\|g_k\| \geq \gamma, \text{ where } \gamma > 0. \tag{19}$$

Eq. (11) can be divided into two parts, as follows:

$$\beta_k^{(1)} = \frac{\|g_k\|^2 - \mu_k |g_k^T g_{k-1}|}{d_{k-1}^T y_{k-1}},$$
$$\beta_k^{(2)} = -\frac{\mu_k}{\alpha_{k-1}} \frac{g_k^T s_{k-1}}{d_{k-1}^T y_{k-1}}.$$

Then the following values can be defined:

$$\xi = \frac{\left\|-g_k + \beta_k^{(2)} d_{k-1}\right\|}{\|d_k\|}, \zeta = \frac{\beta_k^{(1)} \|d_{k-1}\|}{\|d_k\|}.$$

From the definition of $u_k$ it can be derived

$$u_k = \frac{d_k}{\|d_k\|} = \frac{-g_k + (\beta_k^{(1)} + \beta_k^{(2)}) d_{k-1}}{\|d_k\|} = \xi + \zeta \frac{d_{k-1}}{\|d_k\|} = \xi + \zeta u_{k-1}.$$

Since $u_k$ is a unit vector, it follows

$$\|\xi\| = \|u_k - \zeta u_{k-1}\| = \|\zeta u_k - u_{k-1}\|.$$

By using the triangle inequality and $\zeta > 0$, one concludes
$$\|u_k - u_{k-1}\| = 2\|\xi\|. \tag{20}$$

Using the definition of $\xi$, we obtain
$$\|\xi\|\|d_k\| = \left\|-g_k + \beta_{k-1}^{(2)} d_{k-1}\right\| \leq \|g_k\| + \left\|\beta_{k-1}^{(2)}\right\|\|d_{k-1}\|. \tag{21}$$

By using the equations of SWP (Eq. (5)) and line search (Eq. (6)),
$$d_{k-1}^T y_{k-1} \geq (\sigma - 1) g_{k-1}^T d_{k-1}.$$
$$\left|\frac{g_k^T d_{k-1}}{d_{k-1}^T y_{k-1}}\right| \leq \left(\frac{\sigma}{1-\sigma}\right).$$

Thus,
$$\beta_k^{(2)} = -\frac{\mu_k}{\alpha_{k-1}} \frac{g_k^T s_{k-1}}{d_{k-1}^T y_{k-1}} \leq \frac{E}{\alpha_{k-1}} \frac{|g_k^T s_{k-1}|}{d_{k-1}^T y_{k-1}} \leq \frac{E}{\alpha_{k-1}} \frac{\|g_k\|\|s_{k-1}\|}{d_{k-1}^T y_{k-1}}.$$

By using Eq. (21), we obtain the following:
$$\|\xi\|\|d_k\| = \left\|-g_k + \beta_{k-1}^{(2)} d_{k-1}\right\| \leq \|g_k\| + \frac{E}{\alpha_{k-1}} \left|\frac{g_k^T d_{k-1}}{d_{k-1}^T y_{k-1}}\right|\|s_{k-1}\| \leq \gamma + \frac{E}{\alpha_{k-1}}\left(\frac{\sigma}{1-\sigma}\right)B.$$

An application of Eq. (20) leads to
$$\|u_k - u_{k-1}\| = 2\|\xi\| = 2\frac{\gamma + \frac{E}{\alpha_{k-1}}\left(\frac{\sigma}{1-\sigma}\right)B}{\|d_k\|}.$$
$$\|u_k - u_{k-1}\|^2 = 4\frac{\left(\gamma + \frac{E}{\alpha_{k-1}}\left(\frac{\sigma}{1-\sigma}\right)B\right)^2}{\|d_k\|^2}.$$

Utilizing Eq. (19), we obtain the following:
$$\sum_{k\geq 1}^{\infty} \frac{1}{\|d_k\|^2} \leq \infty,$$

which completes the proof.

**Lemma 3.7** Assume that Assumption 1 holds and the sequences $\{g_k\}$ and $\{d_k\}$ are generated using Algorithm 1, where $\alpha_k$ is computed via the WWP line search, such that the sufficient descent condition given in Eq. (7) holds and consider that the method satisfies Property(*). Suppose also that Eq. (19) holds. Then, there exists a constant $\lambda > 0$, such that for any $\Delta \in N$ and any index $k_0$, there exists an index $k > k_0$ that satisfies the following:
$$\left|\kappa_{k,\Delta}^{\lambda}\right| > \frac{\lambda}{2},$$
where $\kappa_{k,\Delta}^{\lambda} = \{i \in N: k \leq i \leq k + \Delta - 1, \|s_i\| > \lambda\}$, $N$ denotes the set of positive integers, and $\left|\kappa_{k,\Delta}^{\lambda}\right|$ denotes the number of elements in $\kappa_{k,\Delta}^{\lambda}$.

From Lemmas 3.5, 3.6, and 3.7, the convergence properties of Algorithm 1 with the SWP line search can be satisfied in a manner similar to that used in Theorem 3.6 presented by Gilbert and Nocedal[7]. Therefore, the proof of the following theorem is omitted.

**Theorem 3.6** Assume that the sequences $\{g_k\}$ and $\{d_k\}$ are generated using Eqs. (1) and (2) with the CG formula $\beta_k^{AZHS}$, and let the step length satisfy Eqs. (4) and (5). If Lemmas 3.5, 3.6, and 3.7 are true, then $\liminf_{k\to\infty} \|g_k\| = 0$.

Note that if Lemma 3.1 does not hold true then it is enough to show that
$$\beta_k^{AZHS} = \frac{\|g_k\|^2 - \mu_k|g_k^T g_{k-1}|}{d_{k-1}^T y_{k-1}},$$
satisfies Property (*) similar to Lemma 3.3 in in [24].

The following theorem shows that if the second case of equation (11) holds i.e.

$$\beta_k^{AZHS} = -\frac{1}{\alpha_k}\mu_k \frac{g_k^T s_{k-1}}{d_{k-1}^T y_{k-1}}, \qquad (22)$$

then we will obtain result in Theorem 3.6.

**Theorem 3.7.** Assume that Assumption 1 holds. Consider the conjugate gradient method in (1) and (2) with Equation (22), where $d_k$ is a descent direction and $\alpha_k$ is obtained by the strong Wolfe line search. Then, the $\liminf_{k\to\infty}\|g_k\|=0$.

**Proof.** We will prove this theorem by contradiction. Suppose Theorem 3.7 is not true. Then equation (19) holds

By squaring both sides of (2), we obtain
$$\|d_k\|^2 = \|g_k\|^2 - 2\beta_k g_k^T d_{k-1} + \beta_k^2 \|d_{k-1}\|^2$$
$$\leq \|g_k\|^2 + 2|\beta_k|\,|g_k^T d_{k-1}| + \beta_k^2 \|d_{k-1}\|^2$$
$$\leq \|g_k\|^2 + \frac{2\mathrm{E}}{\alpha_k}\frac{\|g_k\|\|s_k\|}{(1-\sigma)|g_{k-1}^T d_{k-1}|}(\sigma)|g_{k-1}^T d_{k-1}| + \frac{\mathrm{E}^2}{\alpha_k^2}\frac{\left((\sigma)g_{k-1}^T d_{k-1}\right)^2 \|s_{k-1}\|^2}{(1-\sigma)^2|g_{k-1}^T d_{k-1}|^2}$$
$$\leq \|g_k\|^2 + \frac{2\mathrm{E}}{\alpha_k}\frac{\|g_k\|\|s_k\|}{(1-\sigma)}\sigma + \frac{\mathrm{E}^2}{\alpha_k^2}\frac{(\sigma)^2\|s_{k-1}\|^2}{(1-\sigma)^2}.$$

Further calculation gives
$$\frac{\|d_k\|^2}{\|g_k\|^4} \leq \frac{\|g_k\|^2}{\|g_k\|^4} + \frac{2\mathrm{E}}{\alpha_k}\frac{\|g_k\|\|s_k\|}{(1-\sigma)\|g_k\|^4}\sigma + \frac{\mathrm{E}}{\alpha_k^2}\frac{\sigma^2\|s_{k-1}\|^2}{(1-\sigma)^2\|g_k\|^4}$$
$$\leq \frac{1}{\|g_k\|^2} + \frac{2\mathrm{E}}{\alpha_k}\frac{\|g_k\|\|s_k\|}{(1-\sigma)\|g_k\|^4}\sigma + \frac{\mathrm{E}^2}{\alpha_k^2}\frac{\sigma^2}{(1-\sigma)^2\|g_k\|^4}.$$
$$\leq \frac{1}{\|g_k\|^2} + \frac{2\mathrm{E}}{\alpha_k}\frac{\|s_k\|}{(1-\sigma)\|g_k\|^3}\sigma + \frac{\mathrm{E}^2}{\alpha_k^2}\frac{\sigma^2\|s_{k-1}\|^2}{(1-\sigma)^2\|g_k\|^4}.$$

If
$$\|g_k\|^m = \min\{\|g_k\|^2, \|g_k\|^3, \|g_k\|^4\}, m \in N$$

Then it follows
$$\frac{\|d_k\|^2}{\|g_k\|^4} \leq \frac{1}{\|g_k\|^m}\left(1 + \frac{2\mathrm{E}}{\alpha_k}\frac{\lambda}{(1-\sigma)}\sigma + \frac{\mathrm{E}^2}{\alpha_k^2}\frac{\sigma^2\lambda^2}{(1-\sigma)^2}\right).$$

Also,
$$R = \left(1 + \frac{2\mathrm{E}}{\alpha_k}\lambda\sigma + \frac{\mathrm{E}^2}{\alpha_k^2}\frac{\sigma^2\lambda^2}{(1-\sigma)^2}\right)$$

initiates
$$\frac{\|d_k\|^2}{\|g_k\|^4} \leq \frac{R}{\|g_k\|^m} \leq R\sum_{i=1}^k \frac{1}{\|g_i\|^m}$$
$$\frac{\|g_k\|^4}{\|d_k\|^2} \geq \frac{\varepsilon^m}{kR}.$$

Therefore,
$$\sum_{k=0}^{\infty} \frac{\|g_k\|^4}{\|d_k\|^2} = \infty.$$

This result contradicts (15). Therefore, $\liminf_{k\to\infty}\|g_k\|=0$, completing the proof.

### 4- Numerical results and discussion

To analyse the efficiency of the proposed method, we use more than 200 standard test functions, which are presented in Table 1. The test functions can be obtained from the CUTEst library [29]. The(CUTEr/st) test functions with the SIF extension can be obtained from the website
http://www.cuter.rl.ac.uk/Problems/mastsif.shtml.

The numerical results of CG_Descent5.3 were obtained by running the code provided by Hager and Zhang [30] with memory equal to 0. The numerical results of AZHS are obtained using a modified CG_Descent code via the SWP line search, with $\sigma = 0.1$ and $\delta = 0.01$. If $\mu_k > 1$ then we conclude that $L < 1$ and $\frac{\|g_k\|^2}{|g_k^T g_{k-1}|} > 1$. Thus, it reasonable to modified Eq (11) as follows:

$$\beta_k^{AZHS} = \begin{cases} \frac{\|g_k\|^2 - |g_k^T g_{k-1}|}{d_{k-1}^T y_{k-1}}, & \text{if } \|g_k\|^2 > |g_k^T g_{k-1}| \\ \frac{\|g_k\|^2 - \mu_k |g_k^T g_{k-1}|}{d_{k-1}^T y_{k-1}} - \frac{1}{\alpha_k} \mu_k \frac{g_k^T s_{k-1}}{d_{k-1}^T y_{k-1}}, & \text{if } \|g_k\|^2 > \mu_k |g_k^T g_{k-1}|, \\ -\frac{1}{\alpha_k} \mu_k \frac{g_k^T s_{k-1}}{d_{k-1}^T y_{k-1}}, & \text{otherwise} \end{cases}$$

Note that if $\beta_k^{AZHS} = \frac{\|g_k\|^2 - |g_k^T g_{k-1}|}{d_{k-1}^T y_{k-1}}$, then $\beta_k^{AZHS} \leq \beta_k^{HS}$, thus the proof will be similar to that presented by [7].

The host computer was an AMD A4-7210 APU with AMD Radeon R3 Graphics, installed memory (RAM) of 4 GB, and a 64-bit operating system. The graph of the following results was obtained using SigmaPlot, which is a performance measure introduced by Dolan and Moré [31].

This performance measure was introduced to compare the performance of a set of solvers $S$ on a set of problems $\rho$. For $n_s$ solvers and $n_p$ problems in $S$ and $\rho$, respectively, the measure $t_{p,s}$ is defined as the computation time (e.g., the number of iterations or the CPU time) required for solver $s$ to solve problem $p$.

To create a baseline for comparison, the performance of solver $s$ on problem $p$ is scaled with respect to the best performance of any solver in $S$ on that problem, yielding the ratio

$$r_{p,s} = \frac{t_{p,s}}{\min\{t_{p,s}: s \in S\}}.$$

A parameter $r_M \geq r_{p,s}$ for all $p, s$ is selected such that $r_{p,s} = r_M$ if and only if solver $s$ cannot solve problem $p$. To obtain an overall assessment of the performance of each solver, we define the measure

$$P_s(t) = \frac{1}{n_p} \text{size}\{p \in \rho: r_{p,s} \leq t\}.$$

$P_s(t)$ is the probability for solver $s \in S$ that the performance ratio $r_{p,s}$ will be within a factor $t \in R$ of the best possible ratio. If we define the function $p_s$ as the cumulative distribution function of the performance ratio, then the performance measure $p_s: \mathbb{R} \to [0,1]$ for a given solver is non-decreasing and piecewise continuous from the right. The value of $p_s(1)$ is the probability that the solver will achieve the best performance among all solvers. In general, a solver with higher values of $P_s(t)$, which will lie closer to the upper right corner of the figure, is preferable.

The numerical results are shown in figures 1, 2, 3, and 4. Figure 1 represents the number of iterations showing that the new modification significantly outperforms CG_Descent5.3.

Figure 2 shows that the new modification, AZHS, outperforms CG_Descent5.3 in the number of function evaluations. Figures 3 and 4 represent performance based on the number of gradient evaluations and CPU time, respectively. We can observe that AZHS outperforms CG_Descent 5.3 in CPU time and is significantly competitive with CG_Descent5.3 in the number of function evaluations and gradient evaluations since the later used approximate Wolfe line search with $\sigma = 0.9$ and $\delta = 0.1$. Thus, we can conclude that the $\beta_k^{AZHS}$ outperforms CG_Descent5.3 in all figures.

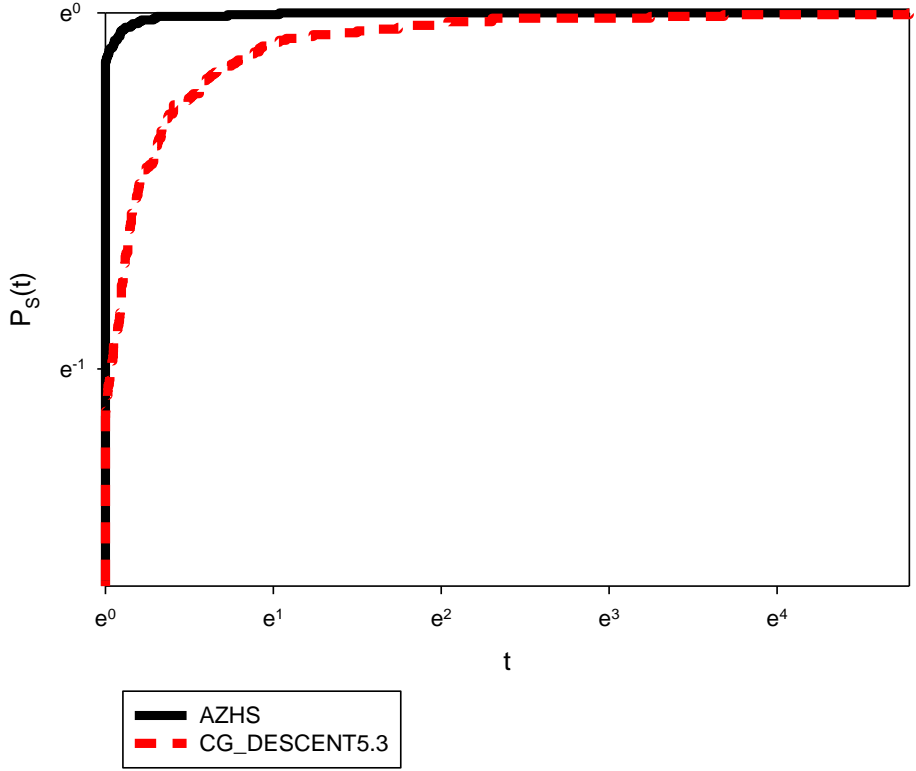

Figure 1. Performance measure based on the number of iterations.

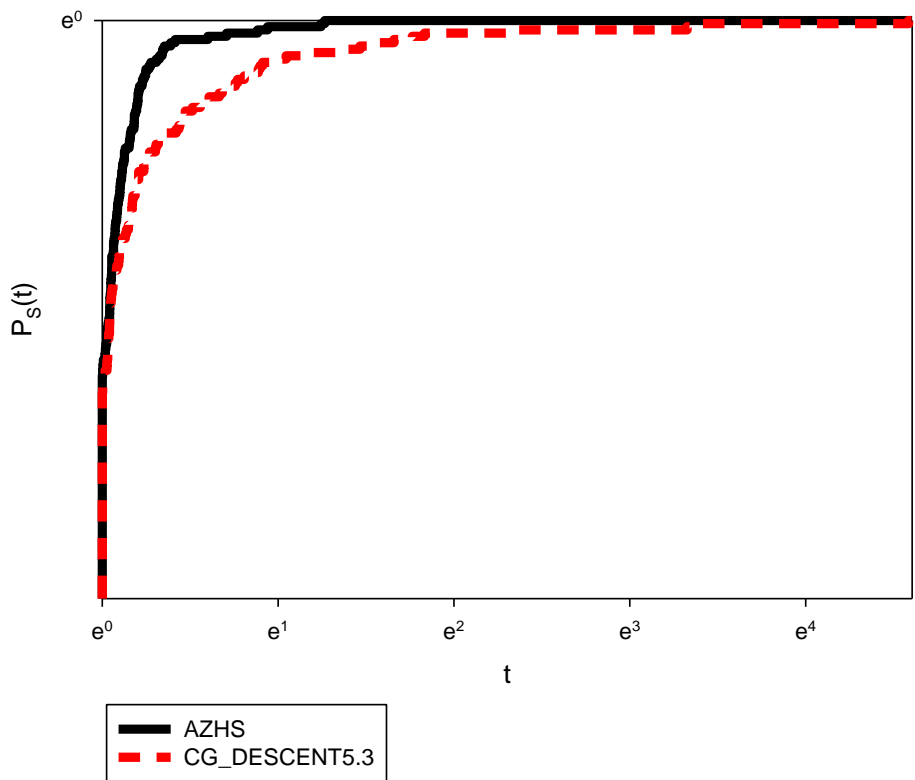

Figure 2. Performance measure based on the function evaluation.

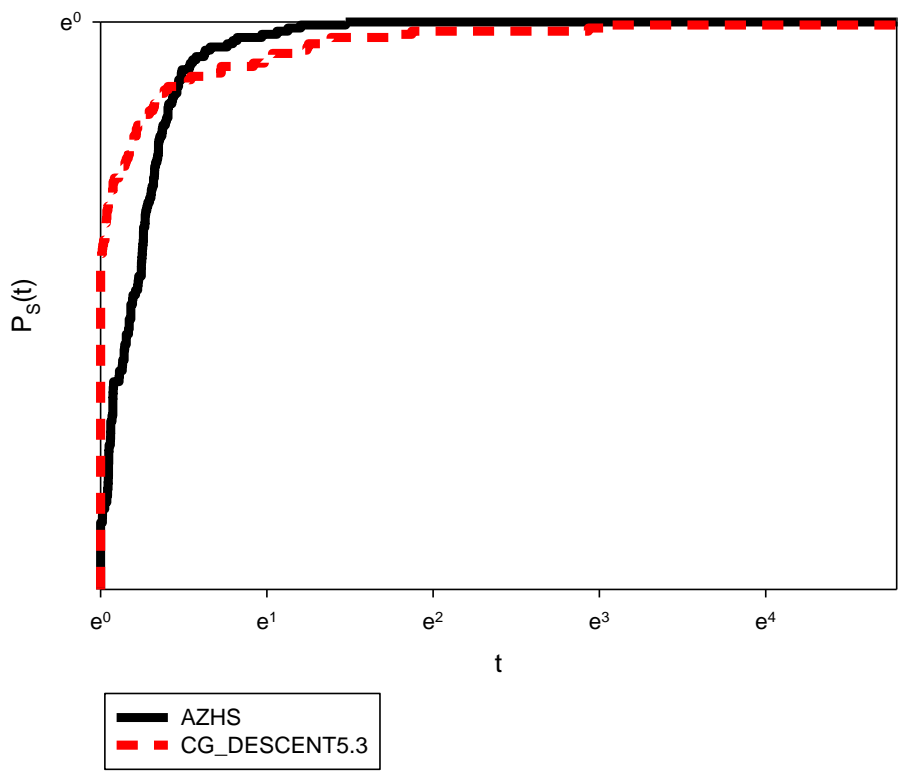

Figure 3. Performance measure based on the gradient evaluation.

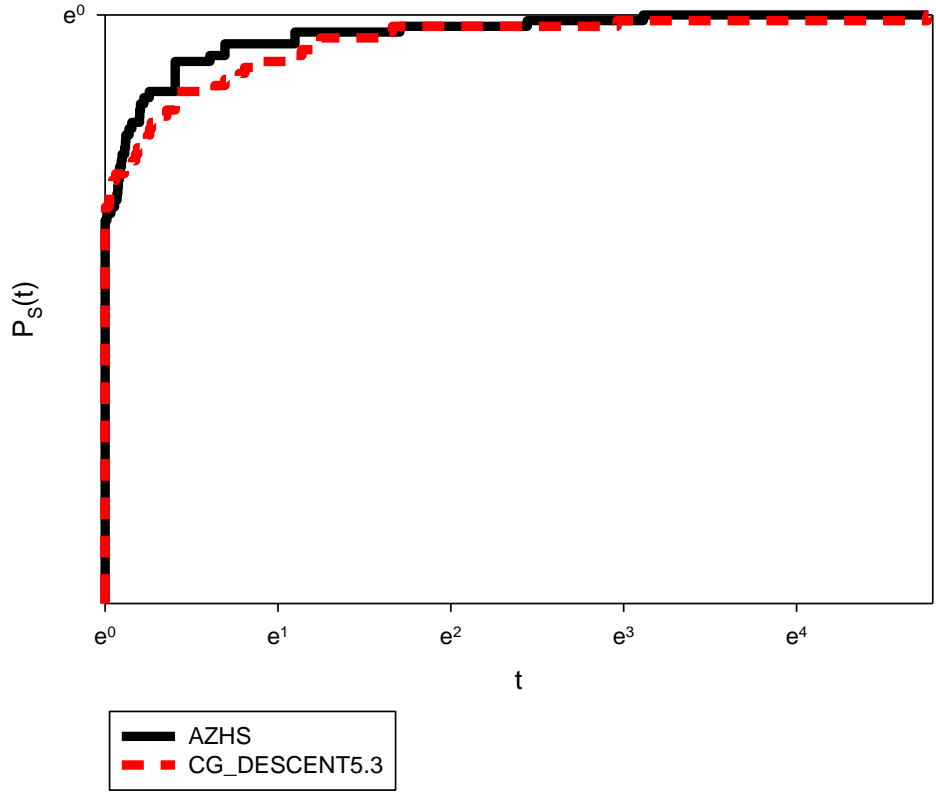

Figure 4. Performance measure based on the CPU time.

Table 1. Test functions

| Function | Dim | Function | Dim | Function | Dim |
|---|---|---|---|---|---|
| AKIVA | 2 | FBRAIN2LS | 4 | OSCIPATH | 10 |
| ALLINITU | 4 | FLETCBV2 | 5000 | PALMER1C | 8 |
| ARGLINB | 200 | FLETCHCR | 1000 | PALMER1D | 7 |
| ARGLINC | 200 | FMINSRF2 | 5625 | PALMER2C | 8 |
| ARWHEAD | 5000 | FMINSURF | 5625 | PALMER3C | 8 |
| BARD | 3 | GENHUMPS | 5000 | PALMER4C | 8 |
| BDEXP | 5000 | GROWTHLS | 3 | PALMER5C | 6 |
| BDQRTIC | 5000 | GULF | 3 | PALMER6C | 8 |
| BEALE | 2 | HAHN1LS | 7 | PALMER7C | 8 |
| BIGGS3 | 6 | HAIRY | 2 | PALMER8C | 8 |
| BIGGS5 | 6 | HATFLDD | 3 | PARKCH | 15 |
| BIGGS6 | 6 | HATFLDE | 3 | PENALTY1 | 1000 |
| BIGGSB1 | 5000 | HATFLDFL | 3 | PENALTY2 | 200 |
| BOX2 | 3 | HATFLDFLS | 3 | PENALTY3 | 200 |
| BOX3 | 3 | HEART6LS | 6 | PENALTY3 | 200 |
| BOX | 10000 | HEART8LS | 8 | POWELLBSLS | 2 |
| BRKMCC | 2 | HELIX | 3 | POWELLSG | 5000 |
| BROYDNBDLS | 10 | HIELOW | 3 | POWER | 10000 |
| BROWNAL | 200 | HILBERTA | 2 | POWERSUM | 4 |
| BROWNBS | 2 | HILBERTB | 10 | PRICE3 | 2 |
| BROWNDEN | 4 | HIMMELBB | 2 | PRICE4 | 2 |
| BROYDN7D | 5000 | HIMMELBF | 4 | QING | 100 |
| BRYBND | 5000 | HIMMELBG | 2 | QUARTC | 5000 |
| CAMEL6 | 2 | HIMMELBH | 2 | RAT43LS | 4 |
| CHNROSNB | 50 | HUMPS | 2 | RECIPELS | 3 |
| CLIFF | 2 | HYDCAR6LS | 29 | ROSENBR | 2 |
| COSINE | 10000 | INDEF | 5000 | ROSENBRTU | 2 |
| CUBE | 2 | INDEFM | 100000 | S308 | 2 |
| CURLY10 | 10000 | INTEQNELS | 12 | SCHMVETT | 5000 |
| CURLY20 | 10000 | JENSMP | 2 | SENSORS | 100 |
| CURLY30 | 10000 | JIMACK | 3549 | SINEVAL | 2 |
| DENSCHNA | 2 | JUDGE | 2 | SINQUAD | 5000 |
| DENSCHNB | 2 | KOWOSB | 4 | SISSER | 2 |
| DENSCHNC | 2 | KSSLS | 1000 | SNAIL | 2 |
| DENSCHND | 3 | LANCZOS1LS | 6 | SPMSRTLS | 4999 |
| DENSCHNE | 3 | LANCZOS2LS | 6 | SROSENBR | 5000 |
| DENSCHNF | 2 | LANCZOS3LS | 6 | SSCOSINE | 5000 |
| DIXMAANA | 3000 | LIARWHD | 5000 | SSI | 3 |
| DIXMAANB | 3000 | LOGHAIRY | 2 | STREG | 4 |
| DIXMAANC | 3000 | LSC1LS | 3 | STRATEC | 10 |
| DIXMAAND | 3000 | LSC2LS | 3 | STRTCHDV | 10 |
| DIXMAANE | 3000 | LUKSAN11LS | 100 | TESTQUAD | 5000 |
| DIXMAANF | 3000 | LUKSAN12LS | 98 | THURBERLS | 7 |
| DIXMAANG | 3000 | LUKSAN13LS | 98 | TOINTGOR | 50 |
| DIXMAANH | 3000 | LUKSAN14LS | 98 | TOINTGSS | 5000 |
| DIXMAANI | 3000 | LUKSAN15LS | 100 | TOINTPSP | 50 |
| DIXMAANJ | 3000 | LUKSAN16LS | 100 | TOINTQOR | 50 |
| DIXMAANK | 3000 | MANCINO | 100 | TQUARTIC | 5000 |
| DIXMAANL | 3000 | MARATOSB | 2 | TRIDIA | 5000 |
| DIXMAANP | 3000 | MEXHAT | 2 | TRIGON1 | 10 |
| DIXON3DQ | 10000 | MEYER3 | 3 | TRIGON2 | 10 |
| DJTL | 2 | MGH09LS | 4 | VANDANMSLS | 22 |
| DMN15332LS | 66 | MGH10LS | 3 | VARDIM | 200 |
| DQDRTIC | 5000 | MGH10SLS | 3 | VAREIGVL | 50 |
| ECKERLE4LS | 3 | MGH17LS | 5 | VESUVIALS | 8 |
| EDENSCH | 2000 | MISRA1BLS | 2 | VESUVIOULS | 8 |
| EGGCRATE | 2 | MISRA1CLS | 2 | VIBRBEAM | 8 |
| EG2 | 1000 | MISRA1DLS | 2 | WAYSEA1 | 2 |
| EIGENALS | 2550 | MODBEALE | 20000 | WAYSEA2 | 2 |
| EIGENBLS | 2550 | MOREBV | 5000 | WOODS | 4000 |
| EIGENCLS | 2652 | MSQRTALS | 1024 | YATP1CLS | 123200 |
| ELATVIDU | 2 | MSQRTBLS | 1024 | YATP2CLS | 123200 |
| ENGVAL1 | 5000 | NCB20 | 5010 | YFITU | 3 |
| ENGVAL2 | 3 | NELSONLS | 3 | ZANGWIL2 | 2 |
| ENSOLS | 9 | NONCVXU2 | 5000 | | |
| EXPFIT | 2 | NONDIA | 5000 | | |
| EXTROSNB | 1000 | NONDQUAR | 5000 | | |
| EXP2 | 2 | OSBORNEA | 5 | | |
| FBRAINLS | 2 | OSBORNEB | 11 | | |

## 5- Application on Heat Conduction Problem [32]

Suppose a rectangular flat plate with dimensions of 5x4 units that generates heat [33]. Suppose the thermal conductivity $k$ is fixed, and the heat production per unit area $f$ is a nonlinear function of the temperature $M$. Our objective is to decide the temperature of the slab such that the temperature beside the perimeter of the slab is zero. Poisson's equation classifies the temperature distribution within this region as follows:

$$k\left[\frac{\partial^2 M}{\partial x^2} + \frac{\partial^2 M}{\partial y^2}\right] + f(M) = 0,$$

If $k = 2$ and $f(M) = 20 - \frac{3}{2}M + \frac{1}{20}M^2$

There are 12 mesh points in total, symmetry reduces the problem to only four distinct temperatures.

$$2(M_2 + M_3 - 4M_1) = -20 + \frac{3}{2}M_1 - \frac{1}{20}M_1^2$$

$$2(M_3 + M_1 + M_4 - 4M_3) = -20 + \frac{3}{2}M_3 - \frac{1}{20}M_3^2$$

$$2(2M_1 + M_4 + 4M_2) = -20 + \frac{3}{2}M_2 - \frac{1}{20}M_2^2$$

$$2(2M_3 + M_2 - 3M_4) = -20 + \frac{3}{2}M_4 - \frac{1}{20}M_4^2$$

These equations, permitted in powers of $M_1$ as follows

$$(M_1^2 - 190M_1) + 40(M_2 + M_3 + 10) = 0,$$

$$M_1 + \frac{M_3^2 - 150M_3 + 400}{40} + M_4 = 0,$$

$$2M_1 + \frac{M_2^2 - 190M_2 + 400}{40} + M_4 = 0,$$

$$(M_4^2 - 150M_4) + 40M_2 + 80M_3 + 400 = 0,$$

The objective function $f$ is constructed by summing the squares of the functions connected with each nonlinear equation as follows:

$$f(M_1, M_2, M_3, M_4, H_1, H_2, H_3, H_4, H_5, H_6) = Q_1 + Q_2 + Q_3 + Q_4,$$

where

$$Q_1 = Q_5^2,$$
$$Q_2 = Q_6^2,$$
$$Q_3 = Q_7^2,$$
$$Q_4 = Q_8^2,$$
$$Q_5 = \frac{1}{20}[M_1^2 + H_1 M_1 + H_2(M_2 + M_3 + H_3)],$$

$$Q_6 = 2[M_1 + \frac{(M_3^2 + H_4 M_3)}{H_2} + H_5 + M_4],$$

$$Q_7 = 2[H_6 M_1 + \frac{(M_2^2 + H_1 M_2)}{H_2} + H_5 + M_4],$$

$$Q_8 = \frac{1}{20}[M_4^2 + H_4 M_4 + H_2 M_2 + H_2 H_6 M_3 + H_2 H_5].$$

If
$$H_1 = -190, H_2 = 40, H_3 = 10, H_4 = -150, H_5 = 10, H_6 = 2,$$
let
$$M_1 = x_1, M_2 = x_2, M_3 = x_3, M_4 = x_4$$

Then, we obtain the following function.

$$f(x_1, x_2, x_3, x_4) = \left(2(x_2 + x_3 - 4x_1) + 20 - 1.5x_1 + \frac{x_1^2}{20}\right)^2 + \left(2(x_1 - 3x_3 + x_4) + 20 - 1.5x_3 + \frac{x_3^2}{20}\right)^2$$

$$+ \left(2(2x_1 + x_4 - 4x_2) + 20 - 1.5x_2 + \frac{x_2^2}{20}\right)^2 + \left(2(x_2 + 2x_3 - 3x_4) + 20 - 1.5x_4 + \frac{x_4^2}{20}\right)^2.$$

We say that $f(x_1, x_2, x_3, x_4)$ Heat Conduction Problem function. By using Algorithm 1 we can find the values of $x_1, x_2, x_3, x_4$ as follows:

$$x_1 = 4.8521,$$
$$x_2 = 6.0545,$$
$$x_3 = 6.4042,$$
$$x_4 = 8.1383.$$

The function value is $1.9631e\text{-}007$

### 6- Application to image restoration

Rrestoring damage images is one of the most important application of CG method. In this study, we applied Gaussian noise with a standard deviation of 25% to the original images in Table 3 after that we use Algorithm 1 to restore these images. To express the efficiency of the proposed method we made a comparison between Algorithm1, CG-Descent5.3 and DL+ in terms of number of iteration, CPU time, and root-mean-square error (RMSE). We utilized RMSE between the restored image and the original true image to calculate the quality of the restored image.

$$RMSE = \frac{\|v - v_k\|_2}{\|v\|}.$$

The restored image is denoted by $v_k$ and the true image by $v$. The RMSE determines the quality of the restored image, in which lower values correspond to higher quality. The results in Table 2 show that the new search direction outperforms CG-Descent5.3 and DL+ in terms of number of iteration, CPU time, and the RMSE value. The criteria for stopping is

$$\frac{\|x_{k+1} - x_k\|_2}{\|x_k\|_2} < \varepsilon.$$

In this context, $\varepsilon = 10^{-3}$. Note that if $\varepsilon = 10^{-4}$ or $\varepsilon = 10^{-6}$, then RMSE remains fixed, meaning that a fixed RMSE can have a variation in the number of iterations.

**Table 2.** Numerical outcomes from images with Gaussian noise with a 25% standard deviation added to the original images using the Dai-Liao CG method, AZHS, as well as CG-Descent5.3.

| Image | Algorithm | Number of iteration | CPU time in seconds | RMSE |
|---|---|---|---|---|
| Mandi 128 pixels | DL+ | 127 | 1.724e+000 | 0.1003 |
| | AZHS | 126 | 1.663e+000 | 0.1002 |
| | CG-Descent5.3 | 134 | 1.825e-001 | 0.1004 |
| Coins 128 pixels | DL+ | 135 | 1.542e+000 | 0.0832 |
| | AZHS | 130 | 1.491e+000 | 0.0824 |
| | CG-Descent5.3 | 133 | 1.491e+000 | 0.0831 |
| Mandi 256 pixels | DL+ | 120 | 1.856e+001 | 0.0519 |
| | AZHS | 111 | 1.545e+001 | 0.0510 |
| | CG-Descent5.3 | 119 | 1.656e+001 | 0.0991 |
| Coins 256 pixels | DL+ | 134 | 1.447e+001 | 0.0506 |
| | AZHS | 120 | 1.164e+001 | 0.0501 |
| | CG-Descent5.3 | 130 | 1.564e+001 | 0.0508 |
| Mandi 512 pixels | DL+ | 114 | 7.981e+001 | 0.0371 |
| | AZHS | 105 | 6.755e+001 | 0.0360 |
| | CG-Descent5.3 | 116 | 7.314e+001 | 0.0472 |
| kids 512 pixels | DL+ | 57 | 6.955e+001 | 0.0377 |
| | AZHS | 56 | 5.325e+001 | 0.0384 |
| | CG-Descent5.3 | 55 | 5.634e+001 | 0.0395 |
| Coins 512 pixels | DL+ | 129 | 7.323e+001 | 0.0326 |
| | AZHS | 128 | 5.248e+001 | 0.0324 |
| | CG-Descent5.3 | 127 | 6.323e+001 | 0.0503 |
| Coins 1024 pixels | DL+ | 128 | 3.441e+002 | 0.0326 |
| | AZHS | 110 | 2.549e+002 | 0.0172 |
| | CG-Descent5.3 | 124 | 2.897e+002 | 0.0289 |

Table 3 shows the outcomes of restoring destroyed images using Algorithm 1, indicating that it can be regarded as an efficient approach.

**Table 3** Restoration of destroyed images of Mandi, Coins, Kides, as well as M.83 by reducing $z$ via Algorithm 1

| Image | Original image | Image with Gaussian noise | Resorted image |
|---|---|---|---|
| Mandi (128 pixels) | 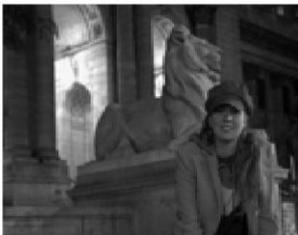 | 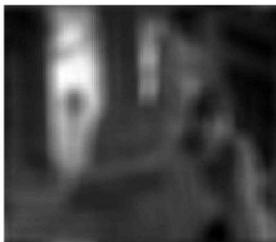 | 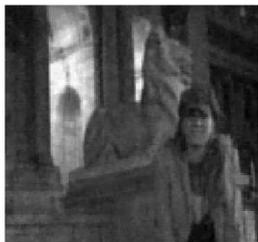 |
| Mandi (256 pixels) | 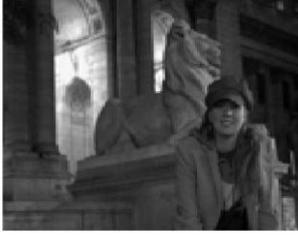 | 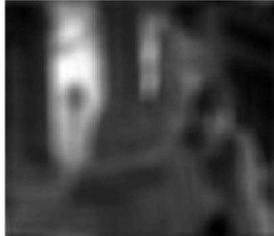 | 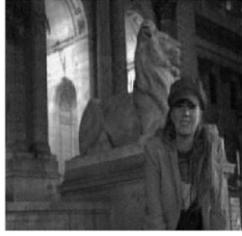 |
| Coins (256 pixels) | 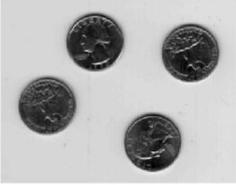 | 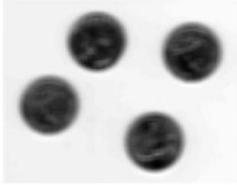 | 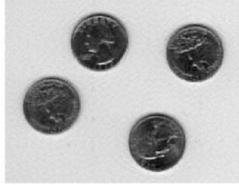 |
| kids(512 pixels) | 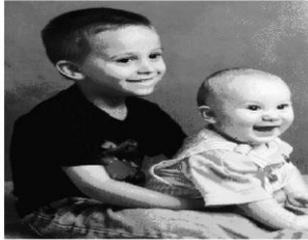 | 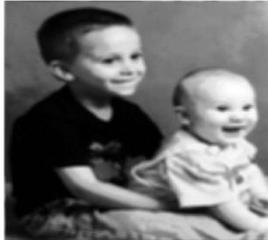 | 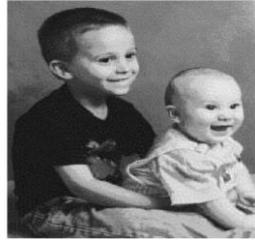 |
| M.83 (1024 pixels) | 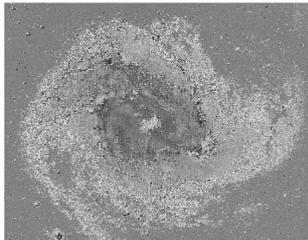 | 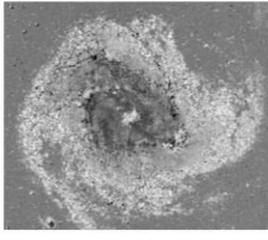 | 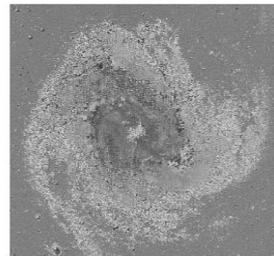 |

### 7- Conclusions

In this study, we investigate a modified HS CG method based on the Dai–Liao conjugacy parameter, with the restart property depending on L. The newly modified CG method inherits global convergence properties and a sufficient descent condition through the SWP line search. Moreover, the numerical results are efficient and competitive with CG Descent5.3. Applications on solving the Heat Conduction Problem and image restoration is presented. In future studies, we will focus on the Lipschitz constant because it plays an essential role in the efficiency and robustness of the CG method.

**Acknowledgements**

This project is supported by the Researchers Supporting Project number (RSP2024R317), King Saud University, Riyadh, Saudi Arabia. We would like to thank Prof. William W. Hager for publishing his code in the implementation of the CG method.

**Availability of data and material**
The data available inside the paper

**Competing interests**
The authors declare that they have no competing interests.

**Authors' contributions**
The authors contributed equally and significantly in writing this paper. All authors read and approved the final manuscript.